\documentclass[12pt,twoside]{article}

\usepackage{amssymb}
\usepackage{amsmath}
\usepackage{bbm}
\usepackage{mathrsfs}
\usepackage{float}
\usepackage{graphicx}
\usepackage{multicol}
\usepackage{pifont}

\makeatletter
\gdef\th@mychange{\normalfont\slshape
   \def\@begintheorem##1##2{\item
        [\hskip\labelsep \theorem@headerfont ##2. ##1  \,--\!--\!--\!--  ]}%
 \def\@opargbegintheorem##1##2##3{%
   \item[\hskip\labelsep \theorem@headerfont ##2. ##1\ {\upshape(}##3{\upshape)}. \,-----  ]}}
\makeatother

\RequirePackage{theorem}
\theoremstyle{mychange}

{\theorembodyfont{\rmfamily}\newtheorem{ttt}{}[section]}
{\theorembodyfont{\rmfamily}}
{\theorembodyfont{\rmfamily}\newtheorem{rem}[ttt]{Remark.}}
{\theorembodyfont{\rmfamily}}
{\theorembodyfont{\rmfamily}}
{\theorembodyfont{\rmfamily}}
{\theorembodyfont{\rmfamily}}
{\theorembodyfont{\rmfamily}}
{\theorembodyfont{\rmfamily}}

{\theorembodyfont{\rmfamily}}

{\theorembodyfont{\itshape}}
{\theorembodyfont{\itshape}}
{\theorembodyfont{\itshape}}
{\theorembodyfont{\itshape}}

{\theorembodyfont{\itshape}}
{\theorembodyfont{\itshape}}
{\theorembodyfont{\itshape}}
{\theorembodyfont{\itshape}}

{\theorembodyfont{\rmfamily}}

{\theorembodyfont{\rmfamily}}

\newcommand{\bbZ}{{\mathbbm Z}}

\newcounter{abc}

\newcounter{iii}

\def\rightend#1#2{{%
 \leavevmode\nobreak\hskip .5em plus 1fil
 \penalty600 \hskip 0pt plus -1filll
 \vadjust{}\nobreak\hskip 0pt plus 1filll%
 #1\parfillskip=#2\relax \par}}

\def\eop{\ifmmode\rule[-22pt]{0pt}{1pt}\ifinner\tag*{$\square$}\else\eqno{\square}\fi\else\rightend{$\square$}{0pt}\fi}

\author{
Andreas-Stephan Elsenhans${}^*$}

\date{}

\title{Rational points on diagonal quartic surfaces}

\begin{document}

\maketitle

\begin{abstract}
We searched up to height $10^7$ for rational points 
on diagonal quartic surfaces.
The computations fill several gaps in earlier lists computed
by Pinch, Swinnerton-Dyer, and Bright.
\end{abstract}

\footnotetext[1]{Mathematisches Institut, Universit\"at Bayreuth, 
Univ'stra\ss e 30, D-95440 Bayreuth, Germany,\\
{\tt Stephan.Elsenhans@uni-bayreuth.de}, 
Website:\! {\tt http://www.staff\!.\!uni-bayreuth.de/$\sim$btm216}\smallskip}

\footnotetext[1]{The author was supported in part by the Deutsche 
Forschungsgemeinschaft (DFG) through a funded research~project.\smallskip}

\footnotetext[1]{${}^\ddagger$The computer part of this work was executed on the servers
of the chair for Computer Algebra at the University of Bayreuth. The author is grateful
to Prof.~M.~Stoll for the permission to use these machines as well as to the system
administrators for their~support.}

\section{Introduction}

The set of rational points on a variety is one of the central objects in
arithmetic geometry. 
For some classes of varieties, one has at least a precise idea how it
should look like.

In the case of Fano varieties, many rational points are expected. 
This expectation is described by the famous conjecture of Manin~\cite{FMT}.
The case of a variety of general type is described by the Lang conjecture.
It claims that the Zariski-closure of the set of rational points has strictly 
smaller dimension than the variety considerd.
Both conjectures are proven only in a few special cases. 
But in the intermediate case
(i.e., varieties with are neither Fano nor of general type) it is not clear how
a general conjecture should look like.

In this note, we inspect diagonal quartic surfaces. These are special K3 surfaces
and form one of the most famous examples for varieties of intermediate type.
More precisely, we focus on surfaces of the form
\begin{eqnarray*}
a x^4 + b y^4 = c z^4 + d w^4 
\end{eqnarray*}
with coefficients $a,b,c,d \in \bbZ$ and $1 \leq a,b,c,d \leq 15$. 
We test local solvability and search for rational points.

\begin{rem}
The only known technique to prove that there are no rational 
points on a K3 surface which has local points everywere 
is given by the Brauer-Manin obstruction. 
The algebraic part of the obstruction was intensively studied by 
Martin Bright in his PhD thesis~\cite{Br1}. 
As explained in~\cite{Br2}, for many diagonal quartic surfaces  
the algebraic and the transcendental Brauer-Manin obstruction can
not explain the absence of rational  points.
\end{rem}

\section{Point search algorithms}
The main idea to search for points of absolute height at most $B$ on 
varieties of the form $f(x,y) = g(z,w)$ is as follows. 

Compute the two sets $\{f(x,y) \mid |x|,|y| \leq B \}$ and
$\{g(z,w) \mid |z|,|w| \leq B \}$. Each solution of the equation
leads to an element in the intersection. On the other hand, one 
can find the solutions when one knows the intersection.
If one can handle sets fast enough, an $O(B^2)$-algorithm results.

But as these sets tend to be very big, a more sophisticated approach has to be used.
One way is given in~\cite{Be} for functions $f$ and $g$ being sums of 
two univariate functions. Then, one can enumerate the two sets above in sorted form. 
From this, the intersection can be formed easily. 

A second approach is
given in~\cite{EJ1} and~\cite{EJ2}. There, the sets are implemented 
using hash tables. To reduce the size of the sets, a page prime $p_p$ is
introduced. Then, the $p_p$ independent problems $f(x,y) = g(z,w)$ and 
$f(x,y) \equiv a\, (\bmod\, p_p)$ are solved for each $a \in [0,p_p-1]$. 
Assuming equidistribution, this reduced the size of the sets approximately
by a factor of $p_p$. 

As we focus on diagonal quartic surfaces, some additional
optimizations can be done. First, one can restrict to non-negative values
of $x,y,z,w$. In about one half of the cases, one can find at least one variable
which must be divisible by $5$.
Further, in many cases, the parity of some or all variables for a primitive
solution can be determined.

Usually, many other modules lead to congruences with could be used for a 
speed up if one knows how to handle them fast on a computer.  
See~\cite{EJ3} for an analysis in a particular example. 

\subsection*{Details of the program written}
The point search was done using the hashing approach. 
In an initialization step, congruences modulo $5$ and powers of $2$ were 
checked to get congruence conditions for primitive solutions. 
The page prime 500083 was chosen and the hash table had 134217728 entries. 
To speed up the modular arithmetic, a table containing fourth roots 
and multiplicative inverse elements modulo the page prime was build 
up in the initialization part. Note that the page prime is congruent to
$3$ modulo $4$ and, thus, the fourth root is unique up to sign.

To avoid multiprecision computations, the computations were done modulo
$2^{64}$. We found less than 100000 coincidences modulo the page prime and 
modulo $2^{64}$. Only these were checked by multiprecision computations.

The running time depends highly on the congruences found. 
Searching on one surface for points up to height $10^7$ took between 
12 and 86 days of CPU time on a 2.27GHz Xeon processor.  
In total, 13 years of CPU time were used.

\section{Results}

In total, there are 7194 quadruples $(a,b,c,d)$ with $a,b,c,d \in \{1,\dots,15\}$,
$a \leq b$, $a \leq c$, and $c \leq d$ and $\gcd(a,b,c,d) = 1$. 
Testing for local solvability excludes 3904 of
the corresponding equations $a x^4 + b y^4 = c z^4 + d w^4$.

A point search with height-bound 10 solves 3009 cases. 
Increasing the bound to 100 leads to solutions for 52 of the 281 remaining equations.
Further, 31 equations have a first solution of height at most 1000.
The remaining 198 equations are the entries of the list~\cite{E2}.
In 21 cases, a solution with height between $10^3$ and $10^4$ was found by Martin Bright.

In 18 cases a new solution with height between $10^4$  and $10^5$ was found.
Further, in 14 cases, there is a first solution with 
height between $10^5$  and $10^6$.
Finally, in 15 cases a first solution with height 
between $10^6$  and $10^7$ was detected.
In 130 cases, still, no solution is known.
A list of them is available at 
{\tt http://www.staff.uni-bayreuth.de/$\sim$btm216/pinch\_list\_cases\_2010.txt}.

As several equations have a first solution with height above
$10^6$, one can not expect the unsolved examples
to be unsolvable. However, $x^4 + y^4 = 6 z^4 + 12 w^4$ is known to be 
unsolvable. See~\cite{Br1} for details. Note that this surface is isomorphic to
$ 3 x^4 + 6 y^4 = 8 z^4 + 8 w^4 $. This is the only pair of 
isomorphic surfaces in the list. 

It would be nice if one could make a complete list of solvable and
unsolvable cases as done in~\cite{CKS} for diagonal cubic surfaces.
\begin{table}[H]
\begin{center}
\begin{tabular}{r|r|r|r||r|r|r|r||r}
$a$ & $b$ & $c$ & $d$ & $x$~~~~ & $y$~~~~ & $z$~~~~ & $w$~~~~ & $\rho$\\  
\hline
1 & 15 &  7 & 11 & 2903019 &  391311 & 1780640 &  549424 & 1 \\
2 & 10 &  7 & 11 & 5742991 & 2277664 & 4262801 & 1865875 & 1 \\
4 & 11 &  7 & 13 &  873483 & 1115876 & 1281143 &  448499 & 1 \\
4 & 12 & 11 & 14 & 3902789 & 1356045 & 1015370 & 2875318 & 1 \\
4 &  5 & 11 & 14 &  394427 & 1355547 & 1112545 &  308333 & 1 \\
5 & 11 &  6 &  7 & 1545359 & 3316097 &  187414 & 3732530 & 1 \\
7 &  9 & 11 & 13 & 3094925 & 7817089 & 6049224 & 6224852 & 1 \\
2 &  9 & 12 & 15 & 3625719 & 1832215 & 1639331 & 2213957 & 1 \\
2 & 11 &  7 &  9 & 2957980 & 1748992 &  468557 & 2308737 & 1 \\  
5 & 14 &  7 &  9 & 1943732 &  493862 &  984595 & 1643257 & 1 \\
4 &  4 & 11 & 13
& 1668661 & 1272265 &  324881 & 1335627 & 1 \\
2 &  3 &  8 & 11 & 1216988 &  924293 &  384555 &  873425 & 1 \\
2 &  8 &  5 & 11 
& 1315404 &  988742 & 1272177 &  470035 & 1 \\
1 &  8 &  4 & 13 & 3730667 & 1735542 & 2189289 & 1913815 & 1  \\
3 &  7 & 12 & 14 
& 1116485 &  269121 &  345539 &  754095 & 2
\end{tabular}
\end{center}
\caption{Surfaces with smallest solutions of height above $10^6$ found.}
\end{table}

\begin{rem}
Some people tend to believe that the arithmetic Picard rank $\rho$ of a K3-surface 
has a strong influence on the set of rational points.
We have no unsolved case with Picard rank greater than 2. The unsolved cases
with rank equal to 2 are 
$[1, 1, 6, 12]$, $[2, 4, 9, 9]$, $[2, 4, 11, 11]$, $[2, 9, 6, 12]$,
$[3, 6, 8, 8]$, $[3, 6, 11, 11]$, $[4, 9, 8, 8]$, and $[6, 12, 11, 11]$. 
On the other hand, the table above contains one equation with Picard rank 2 and
a smallest solution of height $1116485$.

Comparing the rank 1 and the rank 2 cases in the sample, one does not
find a great difference for the proportion of unknown cases.
\end{rem}

\end{document}